\documentclass[12pt,oneside,reqno]{amsart}
\usepackage[latin9]{inputenc}
\usepackage{geometry}
\usepackage{array}
\usepackage{amstext}
\usepackage{amsthm}
\usepackage{amssymb}

\makeatletter

\providecommand{\tabularnewline}{\\}

\numberwithin{equation}{section}
\numberwithin{figure}{section}

\title{}
\def\qed{\quad\vrule height4.17pt width4.17pt depth0pt}

\usepackage{amsfonts}\usepackage{enumerate}

\newtheorem{theorem}{Theorem}\theoremstyle{plain}

\newtheorem{definition}{Definition}\newtheorem{example}{Example}
\newtheorem{lemma}{Lemma}

\newtheorem{proposition}{Proposition}

\numberwithin{equation}{section}

\makeatother

\begin{document}
\title{Enumeration of Bi-Commutative-AG-groupoids}
\author{M. RASHAD$^{A}$}
\author{I. AHMAD$^{A,*}$}
\email{rashad@uom.edu.pk}
\email{iahmaad@hotmail.com}
\address{A. Department of Mathematics, University of Malakand, Chakdara Dir(L),
Pakistan.}
\author{M. SHAH$^{B}$ }
\email{shahmaths\_problem@hotmail.com}
\address{B. Department of Mathematics, Government Post Graduate Collage Mardan,
Pakistan.}
\author{A. Borumand SAEID$^{C}$}
\email{a\_b\_saeid@yahoo.com}
\address{C. Department of Mathematics, Shahid Bahonar University of Kerman
Kerman Iran.}
\keywords{AG-groupoids, AG-test, nuclear square AG-groupoids, transitive groupoids,
alternative AG-groupoids.\\
 {*}Corresponding author}
\begin{abstract}
In this paper we introduce (left, right) bi-commutative AG-group$-$
\\
oids and provide a simple method to test whether an arbitrary AG-groupoid
is bi-commutative AG-groupoid or not. We also explore some of the
general properties of these AG-groupoids. Further we introduce and
study some properties of ideals in these AG-groupoids and decompose
left commutative AG-groupoids by defining some congruences on these
AG-groupoids. 
\end{abstract}

\maketitle

\section{introduction}

An AG-groupoid $S$ is a non-associative algebraic structure that
generalizes the commutative semigroup in general, and satisfies the
left invertive law, $ab\cdot c=cb\cdot a$. One can easily verify
that every AG-groupoid satisfies the medial property, $ab.cd=ac.bd.$
A groupoid $S$ is called a BC Stein groupoid if it satisfies the
BC Stein identity $a\cdot bc=bc\cdot a$ $\forall a,b,c\in S$ \cite{MD,JK}.
We adjoin this property to AG-groupoid to introduce a BC Stein AG-groupoid.
AG-groupoids have been enumerated up to order $6$ in \cite{DSV}.
Using the same techniques and data we also enumerate our BC Stein
AG-groupoids up to order $6.$ Table $1$ provides the enumeration
of our BC Stein AG-groupoids. In this note we provide various non-associative
examples for this class. It is worth mentioning that various other
new classes of AG-groupoids have been recently discovered, enumerated
and discussed by various authors \cite{SAA,SIA,RAAS,SClass,mod,mod m}.
We also define and list various examples of BC Stein AG-groupoid in
Section 2. Section $3$ provides a procedure of testing an arbitrary
AG-groupoid for a BC Stein AG-groupoid. In Section $4$ we discuss
some relations of BC Stein AG-groupoid with already known classes
of AG-groupoids\cite{leftT-2-1,MIM,key-4} and investigate its general
properties. While in Section $5$ we characterize BC Stein AG-groupoids
by the properties of their ideals.

In the following we list some of the already known classes of AG-groupoids
with their identities that will be used frequently in the rest of
this article.

\begin{definition}\label{def 1}\cite{KM}A groupoid $S$ is called
an AG-groupoid if $ab\cdot c=cb\cdot a$ for all $a,b,c\in S.$ \end{definition}

\begin{definition}\label{def 2}\cite{SM,RJT,SAA} An AG-groupoid
$S$ is called --- 

{[}(a){]}
\begin{enumerate}
\item --- a left nuclear square AG-groupoid if $(a^{2}b)c=a^{2}(bc),\,\forall a,b,c\in S.$
\item --- a middle nuclear square AG-groupoid if $(ab^{2})c=a(b^{2}c),\,\,\forall a,b,c\in S.$
\item --- a right nuclear square AG-groupoid if $(ab)c^{2}=a(bc^{2}),\,\,\forall a,b,c\in S.$
\item --- a right alternative AG-groupoid if $a(bb)=(ab)b,\,\,\forall a,b\in S.$
\item --- a Bol$^{*}$-AG-groupoid if $a(bc\cdot d)=(ab\cdot c)d,\,\,\,\forall a,b,c,d\in S.$
\item --- a paramedial AG-groupoid if $ab\cdot cd=db\cdot ca,\,\,\forall a,b,c,d\in S.$ 
\end{enumerate}
\end{definition} \begin{definition}\label{def 3} \cite{MS,RJT,SAA,IMM}An
AG-groupoid $S$ is called --- 

{[}(a){]}
\begin{enumerate}
\item --- an AG$^{*}$-groupoid if $ab\cdot c=b\cdot ac,\,\,\forall a,b,c\in S.$
\item --- a $T_{f}^{4}$-AG-groupoid if $ab=cd\Rightarrow ad=cb,\,\,\forall a,b,c,d\in S.$
\item --- a $T_{b}^{4}$-AG-groupoid if $ab=cd\Rightarrow da=bc,\,\,\forall a,b,c,d\in S.$
\item --- an AG-3-band if $a(aa)=a,\,\,\forall a\in S.$
\item --- a cancellative AG-groupoid if $ax=bx,\textrm{ and }xa=ya\Rightarrow a=b$,
$\forall a,b\in S.$
\item --- a medial AG-groupoid if $ab\cdot cd=ac\cdot bd$, $\forall a,b,c,d\in S.$
\item --- locally associative if $a(aa)=(aa)a,\,\,\forall a\in S.$ 
\end{enumerate}
\end{definition} 
\begin{table}[h]
\begin{centering}
\begin{tabular}{|l|r|r|r|r|}
\hline 
Order & 3 & 4 & 5 & 6\tabularnewline
\hline 
\hline 
Non-associative AG-groupoids & 8 & 269 & 31467 & 40104513\tabularnewline
\hline 
Associative AG-groupoids & 12 & 62 & 446 & 7510\tabularnewline
\hline 
Non-associative BC Stein AG-groupoids & 0 & 0 & 16 & 931\tabularnewline
\hline 
Associative BC Stein AG-groupoids & 5 & 14 & 46 & 173\tabularnewline
\hline 
\end{tabular}
\par\end{centering}
\bigskip{}
 \caption{Enumeration of BC Stein AG-groupoids up to order $6$.}
\end{table}

\medskip{}

\section{BC Stein\textbf{ AG-groupoid}\medskip{}
 }

In this section we define a BC Stein AG-groupoid and give a non-associative
example of lowest order to show its existence.

\noindent \begin{definition}\label{defi 4} An AG-groupoid $S$ is
called a BC Stein AG-groupoid if, 
\begin{eqnarray}
a(bc) & = & (bc)a\,\,\forall a,b,c\in S\label{eq:ST}
\end{eqnarray}
\end{definition} \begin{example} Let $S=\{1,2,3,4,5\}$. Then it
is easy to verify that $S$ is a BC Stein AG-groupoid.\end{example}
\smallskip{}

\begin{center}
\begin{tabular}{c|ccccc}
$\cdot$ & 1 & 2 & 3 & 4 & 5\tabularnewline
\hline 
1 & 1 & 1 & 1 & 1 & 1\tabularnewline
2 & 1 & 1 & 1 & 1 & 1\tabularnewline
3 & 1 & 1 & 1 & 2 & 2\tabularnewline
4 & 1 & 1 & 2 & 2 & 3\tabularnewline
5 & 1 & 1 & 2 & 2 & 3\tabularnewline
\end{tabular}
\par\end{center}

\section{BC Stein\textbf{ AG-test}}

\bigskip{}

In this section we discuss a procedure by \cite{PN} to check an arbitrary
AG-groupoid $(G,\cdot)$ for a BC Stein-AG-groupoid. To this end we
define the following binary operations. 
\begin{eqnarray}
a\circ b & = & a\cdot bx\label{eq: k}\\
a\bigtriangleup b & = & bx\cdot a\label{eq: l}
\end{eqnarray}
The identity $a\cdot bx=bx\cdot a$ holds if, 
\begin{eqnarray}
a\circ b & = & a\bigtriangleup b\label{eq: m}
\end{eqnarray}
To construct table of the operation $``\circ"$ for any fixed $x\in G,$
we multiply turn by turn elements of index column of the $``\cdot"$
table from left with the elements of x-column of the $``\cdot"$ table.
Similarly the table of operation $``\bigtriangleup"$ for any fixed
$x\in G$ is obtained by presenting x-column of the $``\cdot"$ table
as index column and multiplying it with elements of index row of the
$``\cdot"$ table from the left. If the tables for the operation $``\circ"$
and $``\bigtriangleup"$ coincides for all $x\in G,$ then \ref{eq:ST}
holds and the AG-groupoid is a BC Stein AG-groupoid in this case.
We illustrate this procedure in the following example.

\noindent \begin{example} Verify the following AG-groupoid $G$ for
a BC Stein AG-groupoid.\end{example} \smallskip{}

\begin{center}
\begin{tabular}{c|ccccc}
$\cdot$ & 1 & 2 & 3 & 4 & 5\tabularnewline
\hline 
1 & 1 & 1 & 1 & 1 & 1\tabularnewline
2 & 1 & 1 & 1 & 1 & 1\tabularnewline
3 & 1 & 1 & 1 & 1 & 1\tabularnewline
4 & 1 & 1 & 1 & 1 & 2\tabularnewline
5 & 1 & 1 & 4 & 2 & 1\tabularnewline
\end{tabular}\bigskip{}
\par\end{center}

Extend the above table in a way as described above. The tables to
the right of original $``\cdot"$ table are constructed for the operation
$``\circ"$ and the downwards tables are constructed for the operation
$``\triangle".$ It is clear from the extended table that the downward
tables and the tables on the right coincide so $G$ is a BC Stein
AG-groupoid.
\noindent \begin{flushleft}
 
\begin{table}
\noindent \begin{raggedright}
\begin{tabular}{|>{\centering}p{2mm}|>{\centering}p{2mm}>{\centering}p{2mm}>{\centering}p{2mm}>{\centering}p{2mm}>{\centering}p{2mm}|>{\centering}p{2mm}>{\centering}p{2mm}>{\centering}p{2mm}>{\centering}p{2mm}c>{\centering}p{2mm}>{\centering}p{2mm}>{\centering}p{2mm}>{\centering}p{2mm}c>{\centering}p{2mm}>{\centering}p{2mm}>{\centering}p{2mm}>{\centering}p{2mm}c>{\centering}p{2mm}>{\centering}p{2mm}>{\centering}p{2mm}>{\centering}p{2mm}c>{\centering}p{2mm}>{\centering}p{2mm}>{\centering}p{2mm}>{\centering}p{2mm}c}
\hline 
\textbf{\tiny{}$\cdot$} & \textbf{\tiny{}1} & \textbf{\tiny{}2} & \textbf{\tiny{}3} & \textbf{\tiny{}4} & \textbf{\tiny{}5} &  &  &  &  & \multicolumn{1}{c}{} &  &  &  &  & \multicolumn{1}{c}{} &  &  &  &  & \multicolumn{1}{c}{} &  &  &  &  & \multicolumn{1}{c}{} &  &  &  &  & \multicolumn{1}{c|}{}\tabularnewline
\hline 
\textbf{\tiny{}1} & \textbf{\tiny{}1} & \textbf{\tiny{}1} & \textbf{\tiny{}1} & \textbf{\tiny{}1} & \textbf{\tiny{}1} & \textbf{\tiny{}1} & \textbf{\tiny{}1} & \textbf{\tiny{}1} & \textbf{\tiny{}1} & \multicolumn{1}{c}{\textbf{\tiny{}1}} & \textbf{\tiny{}1} & \textbf{\tiny{}1} & \textbf{\tiny{}1} & \textbf{\tiny{}1} & \multicolumn{1}{c}{\textbf{\tiny{}1}} & \textbf{\tiny{}1} & \textbf{\tiny{}1} & \textbf{\tiny{}1} & \textbf{\tiny{}1} & \multicolumn{1}{c}{\textbf{\tiny{}1}} & \textbf{\tiny{}1} & \textbf{\tiny{}1} & \textbf{\tiny{}1} & \textbf{\tiny{}1} & \multicolumn{1}{c}{\textbf{\tiny{}1}} & \textbf{\tiny{}1} & \textbf{\tiny{}1} & \textbf{\tiny{}1} & \textbf{\tiny{}1} & \multicolumn{1}{c|}{\textbf{\tiny{}1}}\tabularnewline
\textbf{\tiny{}2} & \textbf{\tiny{}1} & \textbf{\tiny{}1} & \textbf{\tiny{}1} & \textbf{\tiny{}1} & \textbf{\tiny{}1} & \textbf{\tiny{}1} & \textbf{\tiny{}1} & \textbf{\tiny{}1} & \textbf{\tiny{}1} & \multicolumn{1}{c}{\textbf{\tiny{}1}} & \textbf{\tiny{}1} & \textbf{\tiny{}1} & \textbf{\tiny{}1} & \textbf{\tiny{}1} & \multicolumn{1}{c}{\textbf{\tiny{}1}} & \textbf{\tiny{}1} & \textbf{\tiny{}1} & \textbf{\tiny{}1} & \textbf{\tiny{}1} & \multicolumn{1}{c}{\textbf{\tiny{}1}} & \textbf{\tiny{}1} & \textbf{\tiny{}1} & \textbf{\tiny{}1} & \textbf{\tiny{}1} & \multicolumn{1}{c}{\textbf{\tiny{}1}} & \textbf{\tiny{}1} & \textbf{\tiny{}1} & \textbf{\tiny{}1} & \textbf{\tiny{}1} & \multicolumn{1}{c|}{\textbf{\tiny{}1}}\tabularnewline
\textbf{\tiny{}3} & \textbf{\tiny{}1} & \textbf{\tiny{}1} & \textbf{\tiny{}1} & \textbf{\tiny{}1} & \textbf{\tiny{}1} & \textbf{\tiny{}1} & \textbf{\tiny{}1} & \textbf{\tiny{}1} & \textbf{\tiny{}1} & \multicolumn{1}{c}{\textbf{\tiny{}1}} & \textbf{\tiny{}1} & \textbf{\tiny{}1} & \textbf{\tiny{}1} & \textbf{\tiny{}1} & \multicolumn{1}{c}{\textbf{\tiny{}1}} & \textbf{\tiny{}1} & \textbf{\tiny{}1} & \textbf{\tiny{}1} & \textbf{\tiny{}1} & \multicolumn{1}{c}{\textbf{\tiny{}1}} & \textbf{\tiny{}1} & \textbf{\tiny{}1} & \textbf{\tiny{}1} & \textbf{\tiny{}1} & \multicolumn{1}{c}{\textbf{\tiny{}1}} & \textbf{\tiny{}1} & \textbf{\tiny{}1} & \textbf{\tiny{}1} & \textbf{\tiny{}1} & \multicolumn{1}{c|}{\textbf{\tiny{}1}}\tabularnewline
\textbf{\tiny{}4} & \textbf{\tiny{}1} & \textbf{\tiny{}1} & \textbf{\tiny{}1} & \textbf{\tiny{}1} & \textbf{\tiny{}2} & \textbf{\tiny{}1} & \textbf{\tiny{}1} & \textbf{\tiny{}1} & \textbf{\tiny{}1} & \multicolumn{1}{c}{\textbf{\tiny{}1}} & \textbf{\tiny{}1} & \textbf{\tiny{}1} & \textbf{\tiny{}1} & \textbf{\tiny{}1} & \multicolumn{1}{c}{\textbf{\tiny{}1}} & \textbf{\tiny{}1} & \textbf{\tiny{}1} & \textbf{\tiny{}1} & \textbf{\tiny{}1} & \multicolumn{1}{c}{\textbf{\tiny{}1}} & \textbf{\tiny{}1} & \textbf{\tiny{}1} & \textbf{\tiny{}1} & \textbf{\tiny{}1} & \multicolumn{1}{c}{\textbf{\tiny{}1}} & \textbf{\tiny{}1} & \textbf{\tiny{}1} & \textbf{\tiny{}1} & \textbf{\tiny{}1} & \multicolumn{1}{c|}{\textbf{\tiny{}1}}\tabularnewline
\textbf{\tiny{}5} & \textbf{\tiny{}1} & \textbf{\tiny{}1} & \textbf{\tiny{}4} & \textbf{\tiny{}2} & \textbf{\tiny{}1} & \textbf{\tiny{}1} & \textbf{\tiny{}1} & \textbf{\tiny{}1} & \textbf{\tiny{}1} & \multicolumn{1}{c}{\textbf{\tiny{}1}} & \textbf{\tiny{}1} & \textbf{\tiny{}1} & \textbf{\tiny{}1} & \textbf{\tiny{}1} & \multicolumn{1}{c}{\textbf{\tiny{}1}} & \textbf{\tiny{}1} & \textbf{\tiny{}1} & \textbf{\tiny{}1} & \textbf{\tiny{}1} & \multicolumn{1}{c}{\textbf{\tiny{}2}} & \textbf{\tiny{}1} & \textbf{\tiny{}1} & \textbf{\tiny{}1} & \textbf{\tiny{}1} & \multicolumn{1}{c}{\textbf{\tiny{}1}} & \textbf{\tiny{}1} & \textbf{\tiny{}1} & \textbf{\tiny{}1} & \textbf{\tiny{}1} & \multicolumn{1}{c|}{\textbf{\tiny{}1}}\tabularnewline
\hline 
\textbf{\tiny{}1} & \textbf{\tiny{}1} & \textbf{\tiny{}1} & \textbf{\tiny{}1} & \textbf{\tiny{}1} & \textbf{\tiny{}1} &  &  &  &  &  &  &  &  &  &  &  &  &  &  &  &  &  &  &  &  &  &  &  &  & \tabularnewline
\textbf{\tiny{}1} & \textbf{\tiny{}1} & \textbf{\tiny{}1} & \textbf{\tiny{}1} & \textbf{\tiny{}1} & \textbf{\tiny{}1} &  &  &  &  &  &  &  &  &  &  &  &  &  &  &  &  &  &  &  &  &  &  &  &  & \tabularnewline
\textbf{\tiny{}1} & \textbf{\tiny{}1} & \textbf{\tiny{}1} & \textbf{\tiny{}1} & \textbf{\tiny{}1} & \textbf{\tiny{}1} &  &  &  &  &  &  &  &  &  &  &  &  &  &  &  &  &  &  &  &  &  &  &  &  & \tabularnewline
\textbf{\tiny{}1} & \textbf{\tiny{}1} & \textbf{\tiny{}1} & \textbf{\tiny{}1} & \textbf{\tiny{}1} & \textbf{\tiny{}1} &  &  &  &  &  &  &  &  &  &  &  &  &  &  &  &  &  &  &  &  &  &  &  &  & \tabularnewline
\textbf{\tiny{}1} & \textbf{\tiny{}1} & \textbf{\tiny{}1} & \textbf{\tiny{}1} & \textbf{\tiny{}1} & \textbf{\tiny{}1} &  &  &  &  &  &  &  &  &  &  &  &  &  &  &  &  &  &  &  &  &  &  &  &  & \tabularnewline
\cline{1-6} \cline{2-6} \cline{3-6} \cline{4-6} \cline{5-6} \cline{6-6} 
\textbf{\tiny{}1} & \textbf{\tiny{}1} & \textbf{\tiny{}1} & \textbf{\tiny{}1} & \textbf{\tiny{}1} & \textbf{\tiny{}1} &  &  &  &  &  &  &  &  &  &  &  &  &  &  &  &  &  &  &  &  &  &  &  &  & \tabularnewline
\textbf{\tiny{}1} & \textbf{\tiny{}1} & \textbf{\tiny{}1} & \textbf{\tiny{}1} & \textbf{\tiny{}1} & \textbf{\tiny{}1} &  &  &  &  &  &  &  &  &  &  &  &  &  &  &  &  &  &  &  &  &  &  &  &  & \tabularnewline
\textbf{\tiny{}1} & \textbf{\tiny{}1} & \textbf{\tiny{}1} & \textbf{\tiny{}1} & \textbf{\tiny{}1} & \textbf{\tiny{}1} &  &  &  &  &  &  &  &  &  &  &  &  &  &  &  &  &  &  &  &  &  &  &  &  & \tabularnewline
\textbf{\tiny{}1} & \textbf{\tiny{}1} & \textbf{\tiny{}1} & \textbf{\tiny{}1} & \textbf{\tiny{}1} & \textbf{\tiny{}1} &  &  &  &  &  &  &  &  &  &  &  &  &  &  &  &  &  &  &  &  &  &  &  &  & \tabularnewline
\textbf{\tiny{}1} & \textbf{\tiny{}1} & \textbf{\tiny{}1} & \textbf{\tiny{}1} & \textbf{\tiny{}1} & \textbf{\tiny{}1} &  &  &  &  &  &  &  &  &  &  &  &  &  &  &  &  &  &  &  &  &  &  &  &  & \tabularnewline
\cline{1-6} \cline{2-6} \cline{3-6} \cline{4-6} \cline{5-6} \cline{6-6} 
\textbf{\tiny{}1} & \textbf{\tiny{}1} & \textbf{\tiny{}1} & \textbf{\tiny{}1} & \textbf{\tiny{}1} & \textbf{\tiny{}1} &  &  &  &  &  &  &  &  &  &  &  &  &  &  &  &  &  &  &  &  &  &  &  &  & \tabularnewline
\textbf{\tiny{}1} & \textbf{\tiny{}1} & \textbf{\tiny{}1} & \textbf{\tiny{}1} & \textbf{\tiny{}1} & \textbf{\tiny{}1} &  &  &  &  &  &  &  &  &  &  &  &  &  &  &  &  &  &  &  &  &  &  &  &  & \tabularnewline
\textbf{\tiny{}1} & \textbf{\tiny{}1} & \textbf{\tiny{}1} & \textbf{\tiny{}1} & \textbf{\tiny{}1} & \textbf{\tiny{}1} &  &  &  &  &  &  &  &  &  &  &  &  &  &  &  &  &  &  &  &  &  &  &  &  & \tabularnewline
\textbf{\tiny{}1} & \textbf{\tiny{}1} & \textbf{\tiny{}1} & \textbf{\tiny{}1} & \textbf{\tiny{}1} & \textbf{\tiny{}1} &  &  &  &  &  &  &  &  &  &  &  &  &  &  &  &  &  &  &  &  &  &  &  &  & \tabularnewline
\textbf{\tiny{}4} & \textbf{\tiny{}1} & \textbf{\tiny{}1} & \textbf{\tiny{}1} & \textbf{\tiny{}1} & \textbf{\tiny{}2} &  &  &  &  &  &  &  &  &  &  &  &  &  &  &  &  &  &  &  &  &  &  &  &  & \tabularnewline
\cline{1-6} \cline{2-6} \cline{3-6} \cline{4-6} \cline{5-6} \cline{6-6} 
\textbf{\tiny{}1} & \textbf{\tiny{}1} & \textbf{\tiny{}1} & \textbf{\tiny{}1} & \textbf{\tiny{}1} & \textbf{\tiny{}1} &  &  &  &  &  &  &  &  &  &  &  &  &  &  &  &  &  &  &  &  &  &  &  &  & \tabularnewline
\textbf{\tiny{}1} & \textbf{\tiny{}1} & \textbf{\tiny{}1} & \textbf{\tiny{}1} & \textbf{\tiny{}1} & \textbf{\tiny{}1} &  &  &  &  &  &  &  &  &  &  &  &  &  &  &  &  &  &  &  &  &  &  &  &  & \tabularnewline
\textbf{\tiny{}1} & \textbf{\tiny{}1} & \textbf{\tiny{}1} & \textbf{\tiny{}1} & \textbf{\tiny{}1} & \textbf{\tiny{}1} &  &  &  &  &  &  &  &  &  &  &  &  &  &  &  &  &  &  &  &  &  &  &  &  & \tabularnewline
\textbf{\tiny{}1} & \textbf{\tiny{}1} & \textbf{\tiny{}1} & \textbf{\tiny{}1} & \textbf{\tiny{}1} & \textbf{\tiny{}1} &  &  &  &  &  &  &  &  &  &  &  &  &  &  &  &  &  &  &  &  &  &  &  &  & \tabularnewline
\textbf{\tiny{}2} & \textbf{\tiny{}1} & \textbf{\tiny{}1} & \textbf{\tiny{}1} & \textbf{\tiny{}1} & \textbf{\tiny{}1} &  &  &  &  &  &  &  &  &  &  &  &  &  &  &  &  &  &  &  &  &  &  &  &  & \tabularnewline
\cline{1-6} \cline{2-6} \cline{3-6} \cline{4-6} \cline{5-6} \cline{6-6} 
\textbf{\tiny{}1} & \textbf{\tiny{}1} & \textbf{\tiny{}1} & \textbf{\tiny{}1} & \textbf{\tiny{}1} & \textbf{\tiny{}1} &  &  &  &  &  &  &  &  &  &  &  &  &  &  &  &  &  &  &  &  &  &  &  &  & \tabularnewline
\textbf{\tiny{}1} & \textbf{\tiny{}1} & \textbf{\tiny{}1} & \textbf{\tiny{}1} & \textbf{\tiny{}1} & \textbf{\tiny{}1} &  &  &  &  &  &  &  &  &  &  &  &  &  &  &  &  &  &  &  &  &  &  &  &  & \tabularnewline
\textbf{\tiny{}1} & \textbf{\tiny{}1} & \textbf{\tiny{}1} & \textbf{\tiny{}1} & \textbf{\tiny{}1} & \textbf{\tiny{}1} &  &  &  &  &  &  &  &  &  &  &  &  &  &  &  &  &  &  &  &  &  &  &  &  & \tabularnewline
\textbf{\tiny{}2} & \textbf{\tiny{}1} & \textbf{\tiny{}1} & \textbf{\tiny{}1} & \textbf{\tiny{}1} & \textbf{\tiny{}1} &  &  &  &  &  &  &  &  &  &  &  &  &  &  &  &  &  &  &  &  &  &  &  &  & \tabularnewline
\textbf{\tiny{}1} & \textbf{\tiny{}1} & \textbf{\tiny{}1} & \textbf{\tiny{}1} & \textbf{\tiny{}1} & \textbf{\tiny{}1} &  &  &  &  &  &  &  &  &  &  &  &  &  &  &  &  &  &  &  &  &  &  &  &  & \tabularnewline
\cline{1-6} \cline{2-6} \cline{3-6} \cline{4-6} \cline{5-6} \cline{6-6} 
\end{tabular}
\par\end{raggedright}
{\tiny{}\caption{}
}{\tiny\par}
\end{table}
\par\end{flushleft}

\section{\textbf{Characterization of }BC BC Stein\textbf{ AG-groupoids}\bigskip{}
 }

In this section we find the relations of BC BC Stein AG-groupoids
with already other known classes of AG-groupoids. We start with the
following:

\noindent \begin{lemma}\cite{SM}\label{lem 1}AG$^{**}$-groupoid
is Bol$^{*}$- AG-groupoid.\end{lemma}

\noindent \begin{lemma}\cite{SM}\label{lem 2}Bol$^{*}$-AG-groupoid
is paramedial AG-groupoid.\end{lemma}

\noindent \begin{lemma}\cite{SM}\label{lem 3}AG$^{**}$-groupoid
is left nuclear square AG-groupoid.\end{lemma}

\noindent \begin{proposition}\label{prop 1}Let $S$ be a BC BC Stein
AG-groupoid, then

{[}(a){]} 
\begin{enumerate}
\item $S$ is locally associative AG-groupoid.
\item $S$ is right alternative AG-groupoid.
\item $S$ is AG$^{**}$-groupoid.
\item $S$ is Bol$^{*}$-AG-groupoid.
\item $S$ is paramedial AG-groupoid.
\item $S$ is left nuclear square AG-groupoid.
\item $S$ is right nuclear square AG-groupoid.
\item $S$ is middle nuclear square AG-groupoid.
\item $S$ is nuclear square AG-groupoid.
\end{enumerate}
\end{proposition}

\noindent \emph{Proof.} Let $S$ be a BC BC Stein AG-groupoid, and
$a,b,c\in S.$ 

{[}(a){]}
\begin{enumerate}
\item Since $S$ is a BC BC Stein AG-groupoid, thus $a(bc)=(bc)a$, $\forall a,b,c\in S$.
Thus by replacing $b$ and $c$ by $a$, we get $a(aa)=(aa)a.$ Hence
$S$ is locally associative AG-groupoid.
\item Since $\forall a,b\in S$ we have, $ab.b=(bb)a=a(bb)\Rightarrow(ab)b=a(bb).$
Hence $S$ is right alternative AG-groupoid.
\item Since $S$ is BC BC Stein, thus $\forall a,b,c\in S$, $a(bc)=(bc)a=(ac)b=b(ac)\Rightarrow a(bc)=b(ac).$
Hence $S$ is $AG^{**}$-groupoid.
\item The result follows from Lemma \ref{lem 1}.
\item Follows from Lemma \ref{lem 2}.
\item Follows from Lemma \ref{lem 3}.
\item We show that $S$ is right nuclear square AG-groupoid. Let $a,b,c\in S.$
Then 
\begin{align*}
(ab)c^{2} & =(cc\cdot b)a=(b\cdot cc)a=(bc^{2})a=\\
 & =a(bc^{2})\Rightarrow(ab)c^{2}=a(bc^{2}).
\end{align*}
Hence the result is proved.
\item Let $a,b,c\in S.$ Then 
\begin{align*}
a\cdot b^{2}c & =(bb.c)a=(c\cdot bb)a=(cb^{2})a=\\
 & =ab^{2}\cdot c\Rightarrow a\cdot b^{2}c=ab^{2}\cdot c.
\end{align*}
Hence $S$ is middle nuclear square AG-groupoid.
\item Follows by (\emph{f}), (\emph{g}) and (\emph{h}). \qed
\end{enumerate}
\noindent \begin{theorem} A BC BC Stein AG-groupoid $S$ is a semigroup
if any of the following hold.

{[}(i){]} 
\begin{enumerate}
\item $S$ is AG$^{*}$-groupoid.
\item $S$ is AG-3-band.
\item $S$ is $T^{\,4}$-AG-groupoid.
\item $S$ has a cancellative element. 
\end{enumerate}
\end{theorem}

\noindent \emph{Proof.} Let $S$ be a BC BC Stein AG-groupoid.

{[}(i){]}
\begin{enumerate}
\item Assume that $S$ is AG$^{*}$-groupoid, and $a,b,c\in S.$ Then using
Definitions \ref{def 1}, \ref{def 3} and \ref{defi 4}, we have
\begin{align*}
ab\cdot c & =b.ac=ac.b=bc\cdot a=a\cdot bc\Rightarrow ab\cdot c=a\cdot bc.
\end{align*}
\item Let$S$ be an AG-3-band, and $a,b,c\in S.$ Then by Definitions \ref{def 1},
\ref{def 3} and \ref{defi 4}, we get 
\begin{align*}
ab\cdot c & =cb\cdot a=(cb)(aa\cdot a)=(c\cdot aa)(ba)=\\
 & =(ba)(c\cdot aa)=(bc)(a\cdot aa)=(bc)(aa\cdot a)=\\
 & =bc\cdot a=a\cdot bc\Rightarrow ab\cdot c=a\cdot bc.
\end{align*}
\item Let $S$ be a $T^{4}$-AG-groupoid, and $a,b,c\in S.$ Then by Definitions
\ref{def 1}, \ref{def 3} and \ref{defi 4}, we have 
\begin{align*}
ab\cdot c & =cb.a=a.cb\Rightarrow ab.cb=ac\Rightarrow ac.bb=ac\\
 & \Rightarrow ac\cdot c=a\cdot bb=bb\cdot a=ab\cdot b\Rightarrow ac\cdot b=ab\cdot c\\
 & \Rightarrow bc\cdot a=ab\cdot c\Rightarrow a\cdot bc=ab\cdot c.
\end{align*}
\item Let $S$ has a left cancellative element $x$, and $a,b,c\in S.$
Then by Definitions \ref{def 1}, \ref{def 3} and \ref{defi 4},
it is clear that 
\begin{align*}
x(a\cdot bc) & =(a.bc)x=(x.bc)a=(bc.x)a=(xc\cdot b)a=\\
 & =ab\cdot xc=xc\cdot ab=(ab\cdot c)x=x(ab\cdot c)\\
 & \Rightarrow x(a\cdot bc)=x(ab\cdot c)\Rightarrow a\cdot bc=ab\cdot c.
\end{align*}
\item Let $S$ has a right cancellative element $x$, and $a,b,c\in S.$
Then by Definitions \ref{def 1}, \ref{def 3} and \ref{defi 4},
we have 
\begin{align*}
(a\cdot bc)x & =(x.bc)a=(bc.x)a=ax\cdot bc=ab\cdot xc=xc\cdot ab=\\
 & =(ab\cdot c)x\Rightarrow(a\cdot bc)x=(ab\cdot c)x\Rightarrow a\cdot bc=ab\cdot c.
\end{align*}
Hence $S$ is a semigroup. \qed
\end{enumerate}
\bigskip{}

\section{\textbf{Ideals in }BC BC Stein\textbf{ AG-groupoids}\bigskip{}
 }

A subset $A$ of the AG-groupoid $S$ is a left (resp. right) ideal
of $S$ if, 
\begin{equation}
SA\subseteq A\,(\mbox{\mbox{resp. }}AS\subseteq A)\label{eq:q}
\end{equation}
$A$ is a two sided ideal or simply an ideal of $S$ if it is both
left and right ideal of $S$.

Let $S$ be an AG-groupoid and $A,B\subseteq S$, then $A$ and $B$
are called left connected if $SA\subseteq B$ and $SB\subseteq A.$
Similarly, if $S$ is an AG-groupoid and $A,B\subseteq S$, then $A$
and $B$ are called right connected sets if $AS\subseteq B$ and $BS\subseteq A$
and are called connected if are both left and right connected.

\noindent \begin{lemma} If $S$ is a BC BC Stein AG-groupoid and
$A,B\subseteq S$. Then $AB$ and $BA$ are right and left connected
sets.\end{lemma}

\noindent \begin{proof} Let $S$ be a BC BC Stein AG-groupoid. Then
by Definitions \ref{def 1} and \ref{defi 4}, we get 
\begin{align*}
S(AB) & =(AB)S=(SB)A\subseteq BA\Rightarrow S(AB)\subseteq BA.
\end{align*}
Similarly, 
\begin{align*}
S(BA) & =(BA)S=SA\cdot B\subseteq AB\Rightarrow S(BA)\subseteq AB.
\end{align*}

\noindent Hence $AB$ and $BA$ are left connected sets. Now using
Definition \ref{def 1}, we have 
\begin{align*}
(AB)S=(SB)A\subseteq BA\Rightarrow(AB)\subseteq BA.
\end{align*}
Similarly, 
\begin{align*}
 & (BA)S=SA\cdot B\subseteq AB\Rightarrow S(BA)\subseteq AB.
\end{align*}
Hence $AB$ and $BA$ are right connected sets.\end{proof}

\noindent \begin{proposition}If $S$ is a BC BC Stein AG-groupoid
and $L$ is left ideal of $S$ then $L^{2}$ is an ideal of $S.$\end{proposition}

\noindent \begin{proof}Since $S(L^{2})=(LL)S=(SL)L\subseteq L\cdot L=L^{2}$
and $(L^{2})S=SL\cdot L\subseteq LL=L^{2}.$ Hence $L^{2}$ is an
ideal. \end{proof} \medskip{}

\noindent \begin{theorem} If $S$ is a BC BC Stein AG-groupoid. Then
for any fixed $a$ in $S$

{[}(a){]}
\begin{enumerate}
\item $aS$ is an ideal of $S.$
\item $a(Sa)$ is minimal ideal of $S.$
\item $(aS)a$ is an ideal of $S.$
\end{enumerate}
\end{theorem}

\noindent \emph{Proof.} Let $S$ be a BC BC Stein AG-groupoid. Then

{[}(a){]} 
\begin{enumerate}
\item Let $a$ be any fixed element of $S$ then, 
\begin{eqnarray*}
S(aS) & = & \underset{x,y\in S}{\cup}x(ay)=\underset{x,y\in S}{\cup}(ay)x=\\
 & = & \underset{x,y\in S}{\cup}xy\cdot a=\underset{x,y\in S}{\cup}a\cdot xy\subseteq aS\\
\Rightarrow S(aS) & \subseteq & aS.
\end{eqnarray*}
Thus $aS$ is a left ideal of $S.$ Similarly, 
\begin{align*}
(aS)S & =\underset{x,y\in S}{\cup}(ax)y=\underset{x,y\in S}{\cup}yx\cdot a=\underset{x,y\in S}{\cup}a\cdot yx\subseteq aS\\
\Rightarrow(aS)S & \subseteq aS.
\end{align*}
Thus $aS$ is a right ideal of $S$. Equivalently $aS$ is an ideal.
\item Again by Definitions \ref{def 1}, \ref{defi 4} and Proposition \ref{prop 1}(c),
we get 
\begin{align*}
S(a(Sa)) & =\underset{x,y\in S}{\cup}x(a\cdot ya)=\underset{x,y\in S}{\cup}x(ya\cdot a)=\underset{x,y\in S}{\cup}x(aa\cdot y)=\\
 & =\underset{x,y\in S}{\cup}(aa\cdot y)x=\underset{x,y\in S}{\cup}xy\cdot aa=\\
 & =\underset{x,y\in S}{\cup}a(xy\cdot a)\subseteq a(Sa)\Rightarrow S(a(Sa))\subseteq(a(Sa)).
\end{align*}
Thus $a(Sa)$ is a left ideal of $S.$ Similarly by Definition \ref{def 1},
\ref{defi 4} and Proposition\ref{prop 1}(c), we have 
\begin{align*}
(a(Sa))S & =\underset{x,y\in S}{\cup}(a\cdot xa)y=\underset{x,y\in S}{\cup}(xa\cdot a)y=\\
 & =\underset{x,y\in S}{\cup}(aa\cdot x)y=\underset{x,y\in S}{\cup}yx\cdot aa=\\
 & =\underset{x,y\in S}{\cup}a(yx\cdot a)\subseteq a(Sa)\Rightarrow S(a(Sa))\subseteq(a(Sa))
\end{align*}
Hence $a(Sa)$ right ideal of $S.$ Thus $a(Sa)$ is an ideal of $S$.
Let $L$ be an ideal of $S,$ such that $a\in L.$ Then 
\begin{align*}
a(Sa) & \subseteq L(SL)\subseteq LL\subseteq L\mbox{\mbox{ also} \ensuremath{a(Sa)\subseteq L(SL)=(SL)L\subseteq LL\cdot S\subseteq LS\subseteq L.}}
\end{align*}
Hence$a(Sa)$ is minimal ideal of $S.$
\item Now using Definitions \ref{def 1} and \ref{defi 4}, we have 
\begin{align*}
S((aS)a) & =\underset{x,y\in S}{\cup}x(ay\cdot a)=\underset{x,y\in S}{\cup}x(a\cdot ay)=\underset{x,y\in S}{\cup}(a\cdot ay)x=\\
 & =\underset{x,y\in S}{\cup}(x\cdot ay)x=\underset{x,y\in S}{\cup}(ay\cdot x)a=\underset{x,y\in S}{\cup}(xy\cdot a)a=\\
 & \Rightarrow S((aS)a)=\underset{x,y\in S}{\cup}(a\cdot xy)a\subseteq(aS)a\subseteq(aS)a.
\end{align*}
So $(aS)a$ is left ideal of $S.$ Similarly by Definition \ref{defi 4},
we have 
\begin{align*}
((aS)a)S & =\underset{x,y\in S}{\cup}(ax\cdot a)y=\underset{x,y\in S}{\cup}(a\cdot ax)y=\underset{x,y\in S}{\cup}(y\cdot ax)a=\\
 & =\underset{x,y\in S}{\cup}(ax\cdot y)a=\underset{x,y\in S}{\cup}(yx\cdot a)a=\\
 & =\underset{x,y\in S}{\cup}(a\cdot yx)a\subseteq(aS)a\Rightarrow S((aS)a)\subseteq(aS)a.
\end{align*}
Hence $a(Sa)$ is right ideal. Thus $(aS)a$ is an ideal of $S$.\qed
\end{enumerate}
\medskip{}

\noindent \begin{theorem} Let $L$ be a left ideal and $R$ be a
right ideal of a BC BC Stein-AG-groupoid $S,$ such that $x=x^{2}$
and $y=y^{2}$ for some $x,y\in S.$ Then

{[}(i){]}
\begin{enumerate}
\item $xL=L\cap xS$
\item $Rx=Sx\cap R$
\item $Sy\cap xS=x(Sy)$ 
\end{enumerate}
\end{theorem}

\noindent \emph{Proof. } Let $L$ be a left ideal and $R$ be a right
ideal of $S,$ and $x,y\in S$ such that $x=x^{2}$ and $y=y^{2}$.
Then

{[}(i){]}
\begin{enumerate}
\item Let $t\in xL,$ so $t=xl$ for some $l\in L.$ Since $l\in L$ and
$x\in S$, then $xl\in L,$ so that $t\in L.$ But $l\in L\subseteq S,$
therefore $l\in S,$ and so $xl=t\in xS.$ Hence $t\in L\cap xS,$
that is 
\begin{align}
xL & \subseteq L\cap xS\label{eq:n}
\end{align}
Conversely, let $u\in L\cap xS,$ then $u\in L$ and $u\in xS,$ so
$u=xs$ for some $s\in S.$ Now by Definition \ref{defi 4}, Proposition
\ref{prop 1} (c) and the assumption, we get 
\begin{eqnarray*}
u & = & xs=(xx)s=s(xx)=x(sx)=\\
 & = & x(s\cdot xx)=x(xx\cdot s)=x(xs)=\\
 & = & xu\in xL\Rightarrow u\in xL.
\end{eqnarray*}
Therefore, 
\begin{align}
L\cap xS & \subseteq xL\label{eq:o}
\end{align}
Hence by (\ref{eq:n}) and (\ref{eq:o}), we have $L\cap xS=xL$.
\item Let $b\in Rx,$ so $b=rx$ for some $r\in R.$ Since $r\in R$ and
$x\in S$, then $rx\in R,$ so that $b\in R.$ But $r\in R\subseteq S,$
then $r\in S\Rightarrow rx=b\in Sx,$ hence $b\in R\cap Sx,$ this
implies that 
\begin{align}
Rx & \subseteq R\cap Sx\label{eq:p}
\end{align}
Conversely, let $t\in Sx\cap R,$ then $t\in Sx$ and $t\in R.$ Since
$t\in Sx,$ so $t=sx$ for some $s\in S.$ Now by assumption and Definition
\ref{def 1}, \ref{defi 4}, we have 
\[
t=sx=s(xx)=(xx)s=(sx)x=tx\in Rx\Rightarrow t\in Rx.
\]
Therefore, 
\begin{align}
Sx\cap R & \subseteq Rx\label{eq:q-1}
\end{align}
Hence by equations \ref{eq:p} and \ref{eq:q-1}, we have $Sx\cap R=Rx$.
\item If $a\in Sy\cap xS,$ then for some $s,s'\in S,$ we have $a=sy$
and $a=xs'$. Now by Definitions \ref{def 1}, \ref{defi 4} and assumption,
we have $ay=sy\cdot y=yy\cdot s=s\cdot yy=sy=a.$ and 
\begin{align*}
a & =xs'=s'x\cdot x=s'x\cdot xx=\\
 & =x(s'x\cdot x)=x(xx\cdot s')=x\cdot xs'=xa.
\end{align*}
Therefore $ay=xa=a,$ and so again by Definitions \ref{def 1}, \ref{defi 4}
and assumption, we have 
\begin{align*}
x(ay) & =x(xa)=xx\cdot xa=xa\cdot xx=xx\cdot ax=\\
 & =x\cdot ax=ax\cdot x=xx\cdot a=xa=a.
\end{align*}
$\mbox{\mbox{ Thus \ensuremath{a=x(ay)\Rightarrow}}}a\in x(Sy).$
Therefore,
\begin{center}
\begin{align}
Sy\cap xS & \subseteq x(Sy)\label{eq:r}
\end{align}
\par\end{center}

Conversely, if $a\in x(Sy),$ then by Definitions \ref{def 1}, \ref{defi 4}
and assumption, we have 
\begin{align*}
a & =x(sy)=x(s\cdot yy)=x(yy\cdot s)=\\
 & =x(sy\cdot y)=x(y\cdot sy)=(y\cdot sy)x=(x\cdot sy)y=ay,
\end{align*}
and
\begin{center}
$a=x(sy)=xx\cdot sy=(sy\cdot x)x=(x\cdot sy)x=x(x\cdot sy)=xa.$ 
\par\end{center}

\noindent That is $a\in Sy\cap xS.$ Therefore,
\begin{center}
\begin{align}
x(Sy) & \subseteq Sy\cap xS\label{eq:s}
\end{align}
\par\end{center}

Hence by \ref{eq:r} and \ref{eq:s}, we have $x(Sy)=Sy\cap xS.$
\qed
\end{enumerate}
\medskip{}


\begin{thebibliography}{10}
\bibitem[1]{KM}M.A. Kazim and M. Naseerudin, On almost semigroups,
Portugaliae Mathematica, 36(1) (1977).

\bibitem[2]{MD}M.J. Pelling and D.G. Rogers, Some (2,n) combinatorial
BC BC Stein groupoids, aequationes mathematicae, Volume 21, Issue
1, pp 225-235 (1980).

\bibitem[3]{SM}M. Shah, A theoretical and computational investigation
of AG-groups, PhD thesis, Quaid-i-Azam University Islamabad, (2012).

\bibitem[4]{DSV}A. Distler, M. Shah and V. Sorge, Enumeration of
AG-groupoids, Lecture Notes in Computer Science, Volume 6824/2011,
1-14 (2011).

\bibitem[5]{JK}J. Denes and A.D. Keedwell, Latin Squares and their
applications, Academic Press New York and London, 1974.

\bibitem[6]{MS}Q. Mushtaq and S.M. Yusuf, On Locally Associative
LA-Semigroup, J. Nat. Sci. Math., XIX(1), 57--62 (1979).

\bibitem[8]{RJT}J. R. Cho, J. Jezek and T. Kepka, Paramedial Groupoids,
Czechoslovak Mathematical Journal, 49 (124) 1996.

\bibitem[9]{SAA}M. Shah, I. Ahmad, and A. Ali, On introduction of
new classes of AG-groupoids, Res. J. Recent Sci., Vol. 2(1), 2013,
67-70.

\bibitem[10]{SIA}M. Shah, I. Ahmad, and A. Ali, Discovery of new
classes of AG-groupoids, Res. J. Recent Sci., Vol. 1(11), 2012, 47-49.

\bibitem[11]{RAAS}M. Rashad, I. Ahmad, Amanullah and M. Shah, On
relations between right alternative and nuclear square AG-groupoids,
Int. Mathematica Forum, Vol. 8(5), 237-243, 2013.

\bibitem[12]{PN}P.V. Protic and N. Stevanovic, AG-test and some general
properties of Abel-Grassmann's groupoids, PU. M. A, 6 (1995), 371
- 383.

\bibitem[13]{IMM}I. Ahmad, M. Rashad and M. Shah, Some new result
on T$^{1}$, T$^{2}$ and T$^{4}$-AG-groupoids, Research Journal
of Recent Sciences, Vol. 2(3), 64-66, (2013).

\bibitem[14]{SClass}M. Khan, Faisal, V. Amjad, On some classes of
Abel-Grassmann's groupoids. J. Adv. Res. Pure Math. 3(4), 109--119
(2011)

\bibitem[15]{leftT-2-1} M. Rashad, I. Ahmad and M. Shah, Left transitive
AG-groupoids, www.arXiv.org 1402.5296.

\bibitem[16]{MIM} M. Rashad, I. Ahmad, and M. Shah, Bi-commutative
AG-groupoids, www.arXiv.org.

\bibitem[17]{mod}Amanullah, M. Rashad, I. Ahmad, M. Shah, On Modulo
AG-groupoids, arXiv: 1403.2564

\bibitem[18]{mod m} Amanullah, M. Rashad, I. Ahmad, M. Shah, M. Yousaf,
Modulo Matrix AG-groupoids and Modulo AG-groups arXiv:1403.2304

\bibitem[19]{key-4}P. Yiarayong, On Semiprime and Quasi Semiprime
ideals in AG-groupoids, KKU Res. J. 2013; 18(6): 893-896
\end{thebibliography}
\end{document}